\documentclass[english]{amsart}
\usepackage{amsmath,amssymb,amsthm,amsfonts,amscd,mathrsfs}

\usepackage{babel}
\usepackage[all]{xy}
\usepackage{tikz}
\usetikzlibrary{decorations.pathmorphing}
\usepackage{tikz-cd}
 \usepackage{graphicx}
 \usepackage{epstopdf}
 \usepackage{pb-diagram}
 \usepackage[]{todonotes}
\usepackage{mathtools}

\newcommand{\bdm}{\begin{displaymath}}
\newcommand{\edm}{\end{displaymath}}

\newcommand{\Z}{\mathbb{Z}}

\newcommand{\co}{\colon\thinspace}
\newcommand{\TC}{\mathbf{TC}}

\newcommand{\TCG}{\mathbf{\TC}_G}

\newcommand{\cO}{{\mathcal O}}

\newcommand{\X}{{\mathcal X}}

\newcommand{\A}{{\mathcal A}}
\newcommand{\B}{{\mathcal B}}

\newcommand{\e}{\epsilon}

\newcommand{\orb}{\operatorname{{\it orb}}}

\newcommand{\ev}{\operatorname{{\it ev}}}
\newcommand{\cat}{\operatorname{cat}}

\newcommand{\gcat}{\operatorname{cat_G}}
\newcommand{\secat}{\operatorname{secat}}

\newtheorem{defn}{Definition}[section]

\newtheorem{remark}[defn]{Remark}

\theoremstyle{plain}
\newtheorem{thm}[defn]{Theorem}
\newtheorem{prop}[defn]{Proposition}
\newtheorem{lemma}[defn]{Lemma}
\newtheorem{cor}[defn]{Corollary}

\begin{document}

\title[Morita Invariance]{Morita Invariance of\\
Equivariant Lusternik-Schnirelmann Category \\and Invariant Topological Complexity}
\author{A. Angel}
\address{Departamento de Matem\'aticas y Estad\'istica\\
       Universidad del Norte\\ 
       Km. 5 via Puerto Colombia, Barranquilla, Colombia}    
\address{Departamento de Matem\'aticas\\
          Universidad de los Andes\\
	Carrera 1 N. 18 - 10\\ Bogot\'a, Colombia}
\email{ja.angel908@uniandes.edu.co}

 \author{H. Colman}
 \address{Department of Mathematics, Wright College,
 4300 N.\ Narragansett Avenue, Chicago, IL 60634 USA}
\email{hcolman@ccc.edu}

 \author{M. Grant}
 \address{Institute of Mathematics, University of Aberdeen,
Aberdeen AB24 3UE, United Kingdom}
\email{mark.grant@abdn.ac.uk}
\author{J. Oprea}
\address{ Department of Mathematics, Cleveland State University, Cleveland, Ohio 44115}
\email{j.oprea@csuohio.edu}
\
\date{\today}

\begin{abstract} We use the homotopy invariance of equivariant principal bundles to prove that the equivariant $\A$-category of Clapp and Puppe is invariant under
Morita equivalence. As a corollary, we obtain that both the equivariant Lusternik-Schnirelmann category of a group action  and the  invariant topological
complexity are invariant under Morita equivalence. This allows a definition of topological complexity for orbifolds.
\end{abstract}

\maketitle
\section{Introduction}
Many (and maybe all) orbifolds may be described as global quotients of spaces by compact group actions with finite isotropy groups. But in fact,
the same orbifold may have descriptions involving different spaces and different groups. In order to classify orbifolds, therefore, a notion
of equivalence appropriate to this type of representation must be defined. This notion is provided by \emph{Morita equivalence} (see Section
\ref{sec:morita}). The idea to combine Morita equivalence with standard tools of equivariant algebraic topology (such as Bredon cohomology)
to obtain invariants for orbifolds seems to first appear in \cite{pronk}. In this short paper, we will extend the results of \cite{pronk} to two more
standard equivariant invariants, equivariant Lusternik-Schnirelmann category and ``equivariant'' topological complexity. The word ``equivariant'' is in
quotes here because, in fact, we shall demonstrate that standard equivariant topological complexity is not Morita invariant, so it cannot be a good
invariant of orbifolds. Instead, it is the \emph{invariant} topological complexity that is a Morita invariant. 

Lusternik-Schnirelmann category was
invented around 1930 with the goal of providing a lower bound for the number of critical points for any smooth function on a compact manifold.
For more information, see \cite{Clot03} and Section \ref{sec:LSTC} below. An equivariant version of LS category was given in \cite{Fa, Ma, C}
and the appropriate lower bound for critical points of equivariant functions given. Suffice it to say here that the idea of LS category is to
minimally split a space into atomic elements and to use the number of such atoms as an invariant. Topological complexity was invented in
\cite{Far03} (see Section \ref{sec:LSTC}) with the same idea to split a space into a minimal number of atoms, but with the goal of solving the motion planning problem
in robotics on each atom. Again, the minimal number of atoms is an invariant. Equivariant topological complexity was studied in
\cite{CG} and, in the modified form of invariant topological complexity, in \cite{wacek}. Our main result is that, for compact Lie group actions,
equivariant LS category (even in the extended sense of \cite{Clapp}) and invariant topological complexity are Morita invariant and therefore suitable as
invariants of orbifolds.  This was conjectured in \cite{angelcolman}.

Our method of proof involves the use of equivariant principal bundles and we foresee that these objects may prove
useful in other contexts within the general framework of Morita equivalence of group actions.

\section*{Acknowledgements}
This work was partially supported by grants from the Direcci\'on de Investigaci\'on, Desarrollo e Innovaci\'on, Universidad del Norte (\#2018-18 to Andres Angel) and the Simons Foundation (\#278333 to Hellen Colman and \#244393 to John Oprea). All of the authors wish to thank Wilbur Wright College, Chicago for its hospitality during the Topological Robotics Symposium in February 2018, where much of this work was carried out. The first author acknowledges and thanks the hospitality and financial support provided by the Max Planck Institute for Mathematics in Bonn where part of this paper was  also carried out. 

We would like to thank Michael Farber, Aleksandra Franc, Wac{\l}aw Marzantowicz and Yuli Rudyak for insightful discussions.

\section{Morita equivalence for group actions}\label{sec:morita}
A topological groupoid $G_1\rightrightarrows G_0$ is an internal groupoid in the category of topological spaces; that is a groupoid with a topological space of objects $G_0$ and one of morphisms $G_1$ together with the usual structure maps:  source and target  $s, t\colon  G_1 \to G_0$, identity arrows determined by  $u\colon  G_0 \to G_1$, and composition
$m\colon  G_1 \times_{s, G_0, t} G_1 \to G_1$, all given by continuous maps, such that $s$ (and therefore $t$)  is an open surjection.

Let $G$ be a topological group acting continuously on a Hausdorff space $X$. From this data we can construct a topological groupoid, the {\em translation groupoid} $G\ltimes X$, whose objects are the elements of $X$ and whose morphisms $x \to y$ are pairs $(g, x) \in G \times X$ such that $gx = y$, with composition induced by multiplication in $G$. That is, $G\ltimes X$ is the groupoid $G\times X\rightrightarrows X$ where the source is the second projection and the target map is given by the action.

For each $x\in X$ the {\em isotropy group} $G_x=\{h\in G\mid  hx=x\}$ is a closed subgroup of $G$. The set $Gx=\{gx\mid g\in G\}\subseteq X$ is called the {\em orbit} of $x$ and denoted $[x]_G$.

The translation groupoid $G\ltimes X$ may be regarded as a version of the quotient space $X/G$ which keeps more data than the quotient itself as a topological space.

An {\em equivariant map} $\psi\ltimes \e: G\ltimes X\to H\ltimes Y$ is a pair
$(\psi, \e)$ consisting of a homomorphism $\psi: G \to H$ and a continuous map $\e: X\to Y$,  such that $\e(gx)=\psi(g)\e(x)$ for all $g\in G$, $x\in X$.

\begin{defn} An {\em essential equivalence} $\psi\ltimes \e: G\ltimes X\to H\ltimes Y$ is an equivariant map satisfying:\begin{enumerate}
\item(essentially surjective) $t\circ \pi$ is an open surjection:
$$\xymatrix{X\times_Y(H\times Y) \ar[r]^-{\pi}  \ar[d]^{}& H \times Y \ar[r]^{t}\ar[d]^{p_2}&Y \\ X  \ar[r]^{\e}& Y& }$$
where the square is a pullback.
\item(fully faithful) There is a homeomorphism
$\lambda\colon G \times X \to P$ given by $\lambda(g,x)=(\psi(g),\e(x),x,g x)$ where $P$ is
the following pullback:
$$\xymatrix{P\ar[d]^{}  \ar[r]^{}& H \times Y \ar[d]^{(p_2, t)} \\ X\times X  \ar[r]^{\e\times \e}& Y\times Y }$$

\end{enumerate}
\end{defn}

The first condition implies that for all $y\in Y$, there exists $x\in X$ such that $\e(x)\in [y]_G$,
in other words an essential equivalence is not necessarily surjective but has to reach all of the orbits. The second assures that the isotropies are kept;
an essential equivalence cannot send points from different orbits in $X$ to the same orbit  in $Y$ and there is a bijection between the sets:
$$\{g\in G | x'=gx\}=\{h\in H | \e(x')=h\e(x)\}.$$

\begin{remark}

Let's simplify condition 1). Note that
$$
X\times_Y(H\times Y)  = \{ (x,h,y) \mid \e(x) = s(h,y) \} =\ \{ (x,h,y) \mid \e(x) = y) \}
$$
and we have a homeomorphism
$$
X \times H \rightarrow X \times_Y ( H \times Y)
$$
by sending $(x,h) \rightarrow (x,h,\e(x))$ with inverse given by $(x,h,y) \rightarrow (x,h)$.

Condition 1) then is equivalent to requiring that the map
$$
\varphi : X \times H \rightarrow Y
$$
given by
$(x,h) \rightarrow h \e(x)$ is an open surjective  map.

\end{remark}

\begin{defn} Two actions $G\ltimes X$ and $H\ltimes Y$ are Morita equivalent if there is a third action $K\ltimes Z$ and two essential equivalences
$$\xymatrix{
{G\ltimes X}& {K\ltimes Z} \ar[r]^{(\psi, \e)} \ar[l]_{(\psi', \e')}& {H\ltimes Y}.
}$$
We write $G\ltimes X \sim_M H\ltimes Y$.
\end{defn}

The idea is that two actions are Morita equivalent if they define the same {\em quotient object}, which in general is not just the quotient space.
However, we \emph{do} have the following.

\begin{prop}\label{prop:freeMorita}
Two free actions $G\ltimes X$ and $K\ltimes Y$ are Morita equivalent if and only if their quotient spaces $X/G$ and $Y/K$ are homeomorphic.
\end{prop}

For instance the free actions of $\Z_2$ and $\Z_3$ by rotation on the circle $S^1$ are Morita equivalent but their actions on the disc $D^2$ are not Morita equivalent, $\Z_2\ltimes D^2 \not\sim_M \Z_3\ltimes D^2$, since the isotropy at the origin is not preserved.

By Proposition \ref{prop:freeMorita}, free actions are Morita equivalent if they have the same quotient spaces. The notion of Morita equivalence extends the
idea of ``sameness'' to actions that are not free.

As a topological space, the orbit space of a group action can be very uninteresting, so when talking about the action's orbit space we will  consider extra structure on this
singular space constructed from the translation groupoid. The importance of Morita equivalence lies in the fact that Morita equivalent groupoids yield the same enriched singular space.

Many properties of the action are shared by all Morita equivalent groupoids. Some objects, like orbifolds, are defined as a Morita equivalence class of a certain type of groupoid. In addition, we can often understand a particular invariant of a given group action by analyzing that of a simpler representative of its Morita equivalence class.

In what follows, we prove that the generalized Lusternik-Schnirelmann category of a group action,  $_{\A}cat_G(X)$, is invariant under Morita equivalence. By what we have just said,
for orbifolds defined by a group action, this then provides a well defined invariant as well as a tool to calculate this invariant by reducing the calculations to a simpler group action.

Pronk and Scull  gave a nice characterization of essential equivalences for translation groupoids that can be used in practice to construct or check Morita equivalences. In fact, the argument given for Lie groupoids in \cite{pronk} works for topological groupoids as well. Our main result will extensively use this characterization.

\begin{prop}\cite{pronk}\label{pronk} Any essential equivalence is a composite of maps of the forms (1) and (2) described below.
\begin{enumerate}
\item(quotient map) $G\ltimes X \to G/K\ltimes X/K$ where $K$ is a normal subgroup of $G$ acting freely on $X$.
\item (induction map) $H\ltimes X \to  G\ltimes (G\times_HX)  $ where $H$ is a (not necessarily normal) subgroup of $G$ acting on $X$ and $G\times_HX=G\times X/\sim$ with $(gh^{-1},hx)\sim (g,x)$ for any $h \in H$.
\end{enumerate}
\end{prop}

Notation: if $N$ is a normal subgroup of $X$, we write $N\triangleleft X$.

\begin{remark}\label{rem:quotient}
Note that (1) says that if $G$ acts freely on $X$, then the quotient map $X \to X/G$ is an essential equivalence $G \ltimes X \to \{1\} \ltimes X/G$. This means that
an invariant of a group action is  a Morita invariant only when it is the same for the free $G$-space $X$ and for the quotient $X/G$. Later we will see that this
rules out the standard notion of equivariant topological complexity as a candidate for an analogous orbifold invariant.
\end{remark}

\section{Lusternik-Schnirelmann category and topological complexity}\label{sec:LSTC}

Recall that the {\em Lusternik-Schnirelmann category} of a space $X$, denoted $\cat(X)$, is the least $k$ such that $X$ may be covered  by $k$ open sets
\footnote{Note that some authors only require $k+1$ open sets for $\cat(X)=k$.}
$\{ U_1,\ldots , U_k\}$ such that each inclusion $\iota_{U_i}\co U_i\hookrightarrow X$ is null-homotopic. The open sets $U_i$ are called {\em categorical} sets.

Clapp and Puppe \cite{Clapp} defined a generalization of the LS-category, allowing the open sets to deform into a larger class of objects $\A$, not only into points.

We say that a  subset $U \subseteq X$ is {\em compressible} into a subset $A \subseteq X$, if the inclusion map $\iota_U : U \rightarrow X$ is homotopic to a map $c: U \rightarrow X$ with $c(U) \subseteq A$.

Let $\A$ be a class of spaces, at least one of which is non-empty. The {\em  $\A$-category} of $X$, denoted $_{\A}cat X$, is the least integer $k$ such that $X$ may be covered by $k$  open sets $\{U_1,\ldots ,U_k\}$ such that each inclusion $U_i\hookrightarrow X$ factors up to homotopy through a space $A_i\in \A$. If $\A$ is a class of subsets of $X$, then $_{\A}\cat(X)$ is the least integer $k$ such that $X$ may be covered by $k$ open sets $\{U_1,\ldots ,U_k\}$, each compressible into some $A_i \in \A$.

The {\em sectional category} of a map $p: E\to B$, denoted $\secat(p)$, is the least integer $k$ such that $B$ may be covered by $k$ open sets
$\{ U_1,\ldots, U_k\}$ on each of which there exists a homotopy section of $p$, that is, a map $s: U_i\to E$ such that  $ps\simeq \iota_{U_i}:  {U_i}\hookrightarrow B$.

The sectional category of a fibration $p\co E\to B$  is a lower bound for the category of the base, $\secat(p)\le \cat(B)$ and equality holds if the space $E$ is contractible \cite{S}.

The {\em topological complexity} of a space $X$ is a homotopy invariant defined by  Farber \cite{Far03} in order to study the motion planning problem in robotics. Let $X^I$ denote the space of free paths in $X$, endowed with the compact-open topology. Let $\ev:X^I\to X\times X$ be the evaluation map $\ev(\gamma)=\big(\gamma(0),\gamma(1)\big)$.
A {\em motion planner} on an open subset $U \subseteq X \times X$ is a section of $\ev$ over $U$, i.e. a continuous map $s: {U}\to X^I$ such that  the following diagram commutes:

$$\xymatrix{& X^I \ar[d]^{\ev} \\ {U}  \ar[ur]^{s} \ar[r]^-{\iota_U}& X\times X }$$

\begin{defn} The {\em topological complexity} of a space $X$, denoted $\TC(X)$, is the least integer $k$ such that there exists an open cover of $X\times X$
by $k$ open sets on each of which there is a motion planner.
\end{defn}

\begin{remark}
Note that the homotopy lifting property implies that every homotopy section of a fibration gives an actual section, so commutativity in the diagram
above may be replaced by homotopy commutativity.
\end{remark}

The following theorem provides a characterization of the topological complexity in terms of the sectional category and the ${\A}$-category. Let $\Delta(X)\subseteq X\times X$ denote the diagonal subspace.

 \begin{thm} For a space $X$, the following statements are equivalent:
\begin{enumerate}
\item ${\TC}(X) \leq k$.
\item $\secat(ev) \leq k$: there exist  open sets $U_1,\ldots, U_k$ which cover $X \times X$ and sections $s_i :  U_i \rightarrow X^I$
such that $ev \circ s_i$ is homotopic to the inclusion $U_ i \hookrightarrow X \times X$.
\item $_{\Delta(X)}\cat(X \times X) \leq k$: there exist  open sets $U_1,\ldots, U_k$
which cover $X \times X$ and which are compressible into $\Delta(X)$.
\end{enumerate}
\end{thm}

\section{Equivariant theory}
In this section, we recall some of the equivariant versions of the invariants defined previously.

The {\em equivariant category} of a $G$-space $X$, denoted $\gcat(X)$, is the least integer $k$ such that $X$ may be covered by $k$ invariant open sets $\{ U_1,\ldots , U_k\}$, each of which is $G$-compressible into a single orbit. That is, each inclusion map $\iota_{U_j} : U_j \rightarrow X$ is $G$-homotopic to a $G$-map $c: U_j \rightarrow X$ with $c(U_j) \subseteq [z]_G$ for some $z\in X$.

Let $\A$ be a class of $G$-spaces, at least one of which is non-empty.




The {\em equivariant $\A$-category} of $X$, denoted $_{\A}\cat_G(X)$, is the least integer $k$ such that $X$ may be covered by $k$ invariant open sets $\{U_1,\ldots ,U_k\}$ such that each inclusion $U_i\hookrightarrow X$ factors up to $G$-homotopy through a space $A_i\in \A$. If $\A$ is a class of $G$-invariant subsets of $X$, then $_{\A}\cat_G(X)$ is the least integer $k$ such that $X$ may be covered by $k$ invariant open sets $\{U_1,\ldots ,U_k\}$, each  $G$-compressible into some $A_i \in \A$.
In particular, $_{\A}\cat_G(X)=\cat_G (X)$ when $\A$ is either the class of homogeneous spaces $G/H$, or the class of $G$-orbits of the action.

Note that the $G$-invariance of the objects in $\A$ is implicit in the notation $_{\A}\cat_G(X)$.

The {\em equivariant sectional category} of a $G$-map $p: E\to B$, denoted  $\secat_G(p)$, is the least integer $k$ such that $B$ may be covered by $k$ invariant open sets $\{ U_1,\ldots, U_k\}$ on each of which there exists a $G$-map $s: U_i\to E$ such that  $ps\simeq_G \iota_{U_i}:  {U_i}\hookrightarrow B$.

The first attempt to define an equivariant version of the topological complexity resulted in the following notion. Let  $X$ be a $G$-space. The free path fibration $\ev\co X^I\to X\times X$ is a $G$-fibration with respect to the actions
$$G\times X^I\to X^I, \qquad G\times X\times X\to X\times X,$$
$$g(\gamma)(t)=g(\gamma(t)),\qquad g(x,y)=(gx,gy).$$

\begin{defn}\cite{CG} The {\em equivariant topological complexity} of the
  $G$-space $X$, denoted $\TCG(X)$, is the least integer $k$ such that
  $X\times X$ may be covered by $k$ invariant open sets
  $\{U_1,\ldots ,U_k\}$, on each of which there is a $G$-equivariant
  map $s_i\co U_i\to X^I$ such that  the  diagram commutes:
$$\xymatrix{& X^I \ar[d]^{\ev} \\ {U}_i  \ar[ur]^{s_i} \ar[r]^-{\iota_{U_i}}& X\times X }$$

  \end{defn}
The equivariant topological complexity also has  a characterization in terms of the equivariant sectional and equivariant $\A$-category, just like the non-equivariant one.

\begin{thm} For a $G$-space $X$, the following statements are equivalent:
\begin{enumerate}
\item ${\TC}_G(X) \leq k$.
\item $\secat_G(ev) \leq k$: there exist $G$-invariant open sets $U_1,\ldots, U_k$ which cover $X \times X$ and $G$-equivariant sections $s_i :  U_i \rightarrow X^I$ such that $ev \circ s_i$ is  $G$-homotopic to $U_ i \hookrightarrow X \times X$.
\item {$_{\Delta(X)}\cat_G(X \times X) \leq k$}: there exist $G$-invariant open sets $U_1,\ldots, U_k$ which cover $X \times X$ and are $G$-compressible into $\Delta(X)$.

\end{enumerate}
\end{thm}

Unfortunately, the equivariant topological complexity is \emph{not} invariant under Morita equivalence since there are free actions whose equivariant topological complexity is
different from the topological complexity of the quotient space (see Remark \ref{rem:quotient}). For instance, the equivariant topological complexity of the free
action of $S^1$ on $S^1$ by rotation is at least $2$ since the ordinary topological complexity of the space is a lower bound for the equivariant one,
${\TC}_{S^1}(S^1) \ge {\TC}(S^1)=2$. The quotient space of this action is a point, so ${\TC}(S^1/S^1)=1$.

Another approach to defining a topological complexity in the equivariant setting leads to the definition of the {\em invariant topological complexity}.

Consider the space $X^I \times_{X/G} X^I= \big\{(\alpha,\beta) \in X^I \times X^I : [\alpha(1)]_G=[\beta(0)]_G\big\} $.

The  map $\ev' \colon X^I \times_{X/G} X^I\to X \times X$ given by $ \ev'(\alpha,\beta) = \big(\alpha(0), \beta(1)\big)$ is a $(G\times G)$-fibration with respect to the obvious actions.

\begin{defn}\cite{wacek}
The {\em invariant topological complexity} of $X$, $\TC^G(X)$, is the least integer $k$ such that $X \times X$ may be covered by $k$ $(G \times G)$-invariant open sets
  $\{U_1,\ldots ,U_k\}$, on each of which there is a $(G \times G)$-equivariant section
   $ s_i:{U}_i \to X^I \times_{X/G} X^I$
  such that  the  diagram commutes:
$$\xymatrix{& X^I \times_{X/G} X^I \ar[d]^{\ev'} \\ {U}_i  \ar[ur]^{s_i} \ar[r]^-{\iota_{U_i}}& X\times X }$$
\end{defn}

As in the non-equivariant setting, the commutativity of the diagram can be replaced by the requirement that it commute up to $(G \times G)$-homotopy (see below).

Let $\daleth^{G \times G}(X)$ be the saturation of the diagonal $\Delta(X)$ with respect to the $(G\times G)$-action; that is, $\daleth^{G \times G}(X)$ is
the union of all $(G\times G)$-orbits that meet  $\Delta(X)$.

\begin{thm} For a $G$-space $X$ the following are equivalent:
\begin{enumerate}
\item ${\TC}^G(X) \leq k$.
\item $\secat_{G \times G}(\ev') \leq k$: there exist $(G \times G)$-invariant open sets $U_1,\ldots, U_k$ which cover $X \times X$ and $(G \times G)$-equivariant sections $s_i :  U_i \rightarrow X^I \times_{X/G} X^I$ such that $ev' \circ s_i$ is  $(G \times G)$-homotopic to the inclusion $U_ i \hookrightarrow X \times X$.
\item $_{  \daleth^{G \times G} (X)}\cat_{G\times G}(X \times X) \leq k$: there exist $(G \times G)$-invariant open sets $U_1,\ldots, U_k$ which cover $X \times X$ which are $(G\times G)$-compressible into $\daleth^{G \times G} (X)$.
\end{enumerate}
\end{thm}

We will prove that the invariant topological complexity is invariant under Morita equivalence. Our approach will be to utilize the third characterization of the invariant topological complexity and prove a more general result for the equivariant $\A$-category.

\section{Equivariant $\A$-category and Morita invariance}
From now on, all groups will be compact Lie, and all spaces will be assumed to be metrizable. Under these assumptions, we will prove that the equivariant $\A$-category is invariant under Morita equivalence. Given an equivariant map $\psi\ltimes \e:G\ltimes X \rightarrow H\ltimes Y$ and a $G$-invariant subset $A\subseteq X$, we let $A'\subseteq Y$ denote the saturation of $\e(A)$ with respect to the $H$-action. Then if $\A$ is a class of $G$-invariant subsets of $X$, the class $\A':=\{A' |\;  A\in \A\}$ consists of $H$-invariant subsets of $Y$.

\begin{thm}\label{thm:essential} Let $\psi\ltimes \e:G\ltimes X \rightarrow H\ltimes Y$ be an essential equivalence, and $\A$ a class of $G$-invariant subsets of $X$. Then $$_{\A}\cat_{G}(X) = _{\A'}\cat_{H}(X).$$
\end{thm}

We will prove the result by proving that both maps of type (1) and (2) in Proposition \ref{pronk} preserve the equivariant $\A$-category. Explicitly, we will prove:

\begin{thm}\label{thm:morita}
Let $G$ be a compact Lie group acting on a metrizable space $X$, with subgroup $H\le G$ and normal subgroup $K\lhd G$ such that $K$ acts freely. Let $\A$ be  a class  of $G$-invariant subsets of $X$ and  $\B$ be  a class  of $H$-invariant subsets of $X$. Consider  ${\A}/K=\{A/K \;|\; A\in {\A}\}$ and
$G\times_H {\B}=\{G\times_H B \;|\; B\in {\B}\}$. Then
\medskip
\begin{enumerate}
\item $_{\A}\cat_G X= _{{\A}/K}\cat_{G/K} (X/K)$
\medskip
\item $_{\B}\cat_H X= _{G\times_H {\B}}\cat_{G} (G\times_H X)$.
\end{enumerate}
\end{thm}

Recall that the action of $G/K$ on $X/K$ is given by $(gK)\cdot [x]_K=[gx]_K$ for $g\in G, x\in X$. This is well-defined, since for any $k,k'\in K$ we have
\begin{align*}
(gkK)\cdot [k' x]_K & :=[gkk'x]_K \\
       & = [\ell g x]_K\qquad \mbox{for some $\ell\in K$, by normality}\\
       & = [gx]_K.
\end{align*}
It is then clear that $\A/K$ as defined in Theorem \ref{thm:morita} is a family of $G/K$-invariant subsets of $X/K$. Since the action of $G$ on $G\times_H X$ is given by $g\cdot [g_0,x]=[gg_0,x]$ for all $g,g_0\in G$, $x\in X$, it is also clear that $G\times_H {\B}$ is a family of $G$-invariant subsets of $G\times_H X$.
The proof of Theorem \ref{thm:morita} will follow from Lemmas \ref{lem:quotient1}, \ref{lem:quotient2} and \ref{lem:extension} below.

\begin{lemma}\label{lem:quotient1}
If $K\lhd G$ is a normal subgroup, then $_{\A}\gcat (X) \geq _{\A/K}\cat_{G/K}(X/K)$.
\end{lemma}
\proof

Let $U \subseteq X $ be a $G$-invariant open set $G$-compressible into $A\in \A$ and $H: U\times I \to X $ a $G $-equivariant
homotopy such that $H_0=\iota_U:U\hookrightarrow X$ and $H_1(U)\subseteq A\in \A $.
Consider $V=U/K \subseteq  X /K  $ and define $F: V\times I \to X / K $ as $F([(x)]_{ K }, t)=[H(x,t)]_{ K }$.  The homotopy $F$ is continuous, and:
\begin{enumerate}

\item $F$ is well defined: If $[x]_K=[y]_K$, then there exists $k\in K$ such that $y=kx$.  Since $H$ is $G$-equivariant, then
$F([y]_K, t)=F([kx]_K, t)\stackrel{\text{def}}{=}[H(kx,t)]_K\stackrel{\text{$G$-equiv}}{=}[kH(x,t)]_K\stackrel{\text{$k\in K$}}{=}[H(x,t)]_K\stackrel{\text{def}}{=}F([x]_K, t).$

\item $F$ is $(G/K)$-equivariant: if we denote $\bar g = gK$, then
$F(\bar g([x]_K, t))\stackrel{\text{action}}{=}F([gx]_K, t)\stackrel{\text{def}}{=}
[H(gx,t)]_K\stackrel{\text{$G$-equiv}}{=}[gH(x,t)]_K\stackrel{\text{action}}{=}\bar g[H(x,t)]_K\stackrel{\text{def}}{=}\bar gF([x]_K, t).$

\item $F_0=\iota_U$:

$F_0([x]_K)=[H_0(x)]_K=[x]_K$.

\item $F_1(V)\subseteq A/K$: given $x\in U$,

$F_1([x]_{K})=[H_1(x)]_{K} \subseteq A/K$.
\end{enumerate}
Therefore $V$ is a $(G/K)$-invariant open set $(G/K)$-compressible into $A/K\in \A/K$.

\hfill $\blacksquare$

 In order to prove the inequality $_{\A}\cat_{G}(X) \leq _{\A/K}\cat_{G/K}(X/K)$ we need to recall some of the theory of equivariant principal bundles.
 Here we use the hypothesis that $K$ is a closed (and therefore compact) subgroup of $G$ which acts freely on $X$ (note that freeness of the action was not required in Lemma \ref{lem:quotient1}).

Recall from  \cite[p.56]{tomdieck} the definition of $(\Gamma, \alpha,G)$-bundle: Given $\alpha : \Gamma \rightarrow \operatorname{Aut}(G)$ a continuous group homomorphism between a compact Lie group $\Gamma$ and the automorphism group of a topological group $G$, a $(\Gamma, \alpha,G)$-bundle consists of a locally trivial $G$-principal bundle (with right action!) $p : E \rightarrow B$ together with left $\Gamma$-actions on $E$ and $B$ such that:

\begin{enumerate}
\item  $p$ is $\Gamma$-equivariant;
\item For $\gamma \in \Gamma$, $g\in G$ and $x \in X$ the relation $\gamma(x \cdot g) = (\gamma \cdot x) \alpha(\gamma)(g)$ holds.
\end{enumerate}

The data  $\Gamma, \alpha$, and  $G$  give rise to a semidirect product $\Gamma \times_{\alpha} G$.  The topological space $ \Gamma \times G$  carries the multiplication
$$
(\gamma, g)\left(\gamma^{\prime}, g^{\prime}\right)=\left(\gamma \gamma^{\prime}, \alpha_{\gamma}(g') \cdot g\right)
$$
In particular if $\alpha$ is the trivial homomorphism, we have the group $\Gamma \times G^{op}$ which is isomorphic to $\Gamma \times G$.

The topological group $ \Gamma \times_{\alpha} G$  acts from the left on the total space $E$  of a $(\Gamma, \alpha, G)$-bundle
$$
((\gamma, g), e) \mapsto(\gamma e) g
$$

There is a notion of local triviality for $(\Gamma,\alpha,G)$-bundles which is somewhat complicated to state, so we refer the reader to \cite[p.57]{tomdieck}.
A $(\Gamma,\alpha,G)$-bundle is called \emph{numerable} if it is locally trivial with respect to a numerable cover.

\begin{prop}\label{prop:Kfactor}
Given a  left action of a compact Lie group $G$ on a completely regular space $X$ such that a closed normal subgroup $K \lhd G$ acts freely with $X/K$ paracompact, the quotient map
$p :  X \rightarrow X/K$ is a numerable $(G,\alpha,K)$-bundle, where $\alpha : G \rightarrow \operatorname{Aut}(K)$ is given by conjugation, and $X$ is regarded as a free {\em right} $K$-space via $x \cdot g = g^{-1} x$ for $g\in K$, $x\in X$.
\end{prop}

\proof
Because of freeness of the $K$-action, compactness of $K$ and complete regularity of $X$, we have that $p : X \rightarrow X/K$ is a principal $K$-bundle. See \cite[Proposition II.5.8]{bredon}.

The left $G$-action on $X/K$ is given by $g\cdot [x]_K=[gx]_K$, and so $p$ is obviously $G$-equivariant. For all $\gamma\in G, g\in K$ and $x\in X$ we have
$$
\gamma(x \cdot g) = \gamma ( g^{-1} x ) = \gamma g^{-1} \gamma^{-1} \gamma x = \left ( \gamma g \gamma^{-1} \right ) ^{-1} \gamma x = \left ( \alpha (\gamma)(g) \right )^{-1} \gamma x = (\gamma \cdot x) \alpha(\gamma)(g).
$$
The local triviality follows from \cite[Proposition I.8.10]{tomdieck} by complete regularity  and the numerability from the fact that $X/K$ is paracompact.
\hfill$\blacksquare$






Now let's recall two important structural results about $(\Gamma,\alpha, G)$-bundles. These may be found in section 2 of \cite{MS}.

\begin{prop}[Lemma 2.4 of \cite{MS}]\label{prop:homotopybundles2}
Let $p:E\to B\times I$ be a trivial numerable $(\Gamma,\alpha, G)$-bundle with $\Gamma$ a compact Lie group.
Then there exists a bundle isomorphism $E\to (E | B \times \{0\} ) \times I$ over  $B\times I$ and, hence, $p$ has the
$G$-Covering Homotopy Property.
\end{prop}

\begin{cor}[Lemma 2.5 of \cite{MS}]\label{cor:Gfib}
A numerable $(\Gamma,\alpha, G)$-bundle with $\Gamma$ compact Lie  has the $G$-Covering Homotopy Property(i.e. is a $G$-fibration).
\end{cor}

We are now ready to prove:

\begin{lemma}\label{lem:quotient2}
Suppose that $G$ is a compact Lie group acting on a metrizable space $X$.
If $K\lhd G$ is a normal closed subgroup which acts freely, then $_{\A}\cat_{G}(X) \leq _{\A/K}\cat_{G/K}(X/K)$.
\end{lemma}

\proof

Recall that if $X$ is a metrizable space then it is completely regular and paracompact. If a compact Lie group $G$ acts on a metrizable $X$ then $X /G$ is metrizable \cite[Proposition 1.1.12]{palais} and therefore paracompact.

Note that, by Proposition \ref{prop:Kfactor},  $p :  X \rightarrow X/K$ is a numerable $(G,\alpha,K)$-bundle and so is a $G$-fibration by
Corollary \ref{cor:Gfib}.
Let $V\subseteq X/K$ be a $(G/K)$-invariant open set compressible into $A/K\in \A/K$
and $F: V\times I \to X/K$ a $(G/K)$-equivariant homotopy such that $F_0=\iota_V:V\hookrightarrow X/K$ and $F_1(V)\subseteq A/K\in \A/K$.

Consider $U=p^{-1}(V)$ where $p: X\to X/K$ is the projection onto the quotient space.
Because the $G$-action on $V$ factors through the $G/K$-action, the map $F$ is also a $G$-equivariant map.  Since $p$ is a $G$-fibration,
we can find a $G$-lift $\widetilde F$ of the composition $U \times I \to V \times I \to X/K$,
$$\xymatrix{
U \times 0 \ar[r]^-{\mathrm{incl}} \ar[d]_-{p \times 0} & X \ar[d]_-p \\
U \times I \ar[ur]^-{\widetilde F} \ar[r]^-{F (p\times 1)} & X/K.
}
$$
Because $F_1(V) \subseteq A/K$, we see that $\widetilde F_1(U) \subseteq A$ (since $\A$ consists of $G$-invariant sets). Hence the number of sets in
an $\A/K$-categorical covering of $X/K$ gives an upper bound for the minimal number of $\A$-categorical sets covering $X$.




\hfill$\blacksquare$

We have therefore proved Theorem \ref{thm:morita} (1).  Part (2) is proved by applying part (1) twice.

\begin{lemma}\label{lem:extension} If $G$ is a compact Lie group, $H\le G$ is a closed subgroup of $G$ acting on a metrizable space $X$, then
$$_{\B}\cat_H X=_{G \times_H \B}\gcat (G\times_H X).$$
\end{lemma}
\proof
Consider the action of the compact Lie group $G \times H$ on $G \times X$ given by
$$
(g,h) \cdot (\overline{g},x) = (g \overline{g} h^{-1},hx).
$$
Note that $G \times X$ is metrizable because it is the product of metrizable spaces and
$G$ is a closed normal subgroup of  $G \times H$ such that the action of  $G$ on $G \times X$  is free. The quotient $(G \times X )/G$ is homeomorphic to $X$.

Similarly  $H$ is a closed normal subgroup of $ G \times H$ and the action of $H$ on $G \times X$ is free with quotient $(G \times X ) /H = G \times_H X$.

We will apply theorem \ref{thm:morita} (1.) with a class of $G \times H$-invariant subsets coming from the $H$-invariant collection $\B$.

Given a collection  $\B$ of $H$-equivariant subsets of $X$, let $G \times \B = \{ G \times B \mid B \in \B\}$  be the corresponding collection
of $G \times H$-equivariant subsets.  Now let's apply twice what we have proved above. We have the following:

\begin{enumerate}
\item
$_{ \left ( G \times \B \right ) /H}\cat_{ \left ( G \times H \right )  /H}( \left ( G \times X \right )  / H ) = _{  G \times \B }\cat_{G \times H} (G \times X )$
\item $ _{  G \times \B }\cat_{G \times H} (G \times X )= _{ \left ( G \times \B \right ) /G}\cat_{ \left ( G \times H \right )  /G}(\left ( G \times X \right ) / G )$.
\end{enumerate}
But we also have
\begin{enumerate}
\item $ \left ( G \times \B \right )/ G= \B$
\medskip
\item $\left ( G \times \B \right )  / H = G \times_H \B$,
\end{enumerate}
and therefore
\begin{align*}
_{ \left ( G \times_H \B \right )  }\cat_{ G }(  G \times_H X  ) & = _{ \left ( G \times \B \right ) /H}\cat_{ \left ( G \times H \right )  /H}( \left ( G \times X \right )  / H ) \\
& = _{  G \times \B }\cat_{G \times H} (G \times X ) \\
& = _{ \left ( G \times \B \right ) /G}\cat_{ \left ( G \times H \right )  /G}(\left ( G \times X \right ) / G ) \\
 & = _{ \B }\cat_{  H }( X  ).
\end{align*}

\hfill $\blacksquare$

Many $G$-spaces are of the form $G \times_H X$ for $X$ an $H$-space.  Consider a $G$-space $Y$ with a $G$-equivariant map $f: Y \rightarrow G/H$.
If we take $X=f^{-1}(eH)$, then the natural map
$$
 F : G \times_H X \rightarrow Y
$$
is a $G$-equivariant homeomorphism (for $G$ a compact Lie group, $H$ closed) and we have
$$
\cat_G (X) = \cat_H ( f^{-1}(eH) )
$$

\noindent
 In fact the map $F$  gives rise to an equivalence of categories
$$
G-Top_{G/H} \sim H-Top
$$
between  the category of $G$-spaces over $G/H$ and the category $H$-spaces \cite[Page 33]{tomdieck}. Under this equivalence the generalized $G$-equivariant
category of metrizable $G$-spaces over $G/H$ corresponds to the generalized $H$-equivariant category of metrizable $H$-spaces.

\begin{cor}
The equivariant category is invariant under Morita equivalence for compact Lie group actions on metrizable spaces.
\end{cor}
\proof
Let $\psi\ltimes \e: G\ltimes X\to H\ltimes Y$ be an essential equivalence. We have that
$_{\A}\cat_G(X)=\cat_G (X)$ when $\A$ is the class of $G$-orbits. In this case, the class $\A'=\{ A' |\;  A\in \A\}$ is the class of orbits of the action of $H$ on $Y$ since $\psi\ltimes \e$ is an essential equivalence and the result follows from Theorem \ref{thm:essential} above.\hfill$\blacksquare$

\begin{cor}\label{MoritaInvariant}
The invariant topological complexity is invariant under Morita equivalence for compact Lie group actions on metrizable spaces.
\end{cor}
\proof
If $\psi\ltimes \e: G\ltimes X\to H\ltimes Y$ is an essential equivalence, then
$(\psi\times \psi)\ltimes (\e\times \e): (G\times G)\ltimes (X\times X)\to (H\times H)\ltimes (Y\times Y)$ is an essential equivalence.
Since $\daleth^{H \times H}(Y)$ is the saturation with respect to the $(H\times H)$-action
 of $\e(\daleth^{G \times G}(X))$, we have that
 $${\TC}^G(X)= _{  \daleth^{G \times G} (X)}\cat_{G\times G}(X \times X)=_{  \daleth^{H \times H} (Y)}\cat_{H\times H}(Y \times Y)={\TC}^H(Y).$$
 \hfill$\blacksquare$

\begin{remark}
It is likely that everything we have done above for compact Lie group actions also holds for proper actions of discrete groups on a restricted
class of spaces (such as ANRs). Because our main focus is on developing an invariant for orbifolds, we have however chosen to concentrate on the
case of compact Lie groups.
\end{remark}

Higher topological complexities were introduced by Rudyak \cite{rudyak}. In \cite{bayehsarkar} higher analogs of equivariant and invariant topological complexity are introduced.  Let $\Delta_n(X) \subseteq X^n$ be the diagonal and let $\daleth_n(X)$ be the saturation of $\Delta_n(X)$ with respect to the $G^n$-action on $X^n$. Define the
$n^{\mathrm{th}}$-\emph{higher equivariant topological complexity} of $X$ by
$$
\TC_{G,n} = _{  \Delta_n (X)}\cat_{G^n}(X^n)
$$
and the $n^{\mathrm{th}}$-\emph{higher invariant} topological complexity of $X$ by
$$
\TC^{G,n} = _{  \daleth_n (X)}\cat_{G^n}(X^n).
$$
The higher equivariant topological complexities are not Morita invariant. For instance, the higher equivariant topological complexity of the free action of $S^1$ on $S^1$ by rotation is at least $2$ since the ordinary topological complexity of the space is a lower bound for the equivariant ones
$$
2=\TC(S^1) \leq \TC_n(S^1)\leq \TC_{S^1,n}(S^1).
$$
The quotient space of this action is a point, so ${\TC_{e,n}}(S^1/S^1)=1$.

\begin{cor}
The higher invariant topological complexities are invariant under Morita equivalence for compact Lie group actions on metrizable spaces.
\end{cor}
\proof
If $\psi\ltimes \e: G\ltimes X\to H\ltimes Y$ is an essential equivalence, then
$\prod \psi \ltimes \prod \e: G^n\ltimes X^n \to H^n\ltimes Y^n$ is an essential equivalence.

Since $\daleth_{n}(Y)$ is the saturation with respect to the $(H^n)$-action
 of $\e(\daleth_{n}(X))$, we have that
 $${\TC}^{G,n}(X)= _{  \daleth_{n} (X)}\cat_{G^n}(X^n)=_{  \daleth_{n} (Y^n)}\cat_{H^n}(Y^n)={\TC}^{H,n}(Y).$$

 \hfill$\blacksquare$

Note that $\cat(X^G) \leq \gcat(X)$ and we have
$$\max\{\cat(X^G), \cat(X/G)\} \leq \gcat(X).$$
But since equivariant LS-category is Morita invariant and $\cat(X^G)$ is not, by changing by a Morita invariant action we can improve the lower bound.
In fact, by Lemma \ref{lem:quotient1}, for $K \lhd G$ (not necessarily acting freely) we have
$$
\cat_{G/K}(X/K) \leq \gcat (X)
$$
and by taking fixed points
$$
\cat( (X/K)^{G/K} ) \leq \cat_{G/K}(X/K).
$$
This then gives the inequality
$$
\max_{ K \lhd G} \{\cat( (X/K)^{G/K})\} \leq \gcat(X)
$$
which generalizes the known inequality (taking  $K=\{e\}$ and $K=G)$.

There are similar inequalities for invariant TC and the higher analogs. In \cite[ Remark 3.9, Corollary 3.26]{wacek} it is shown that
$$\max\{\TC(X^G), \TC(X/G)\} \leq \TC^G(X),$$
and in \cite[proposition 4.4 and proposition 4.9]{bayehsarkar} that
$$\max\{\TC_n(X^G), \TC_n(X/G)\} \leq \TC^{G,n}(X).$$

\noindent By the same argument as before,
$$
\max_{ K \lhd G} {\TC( (X/K)^{G/K})} \leq \TC^G(X)
$$
and
$$
\max_{ K \lhd G} {\TC_n( (X/K)^{G/K})} \leq \TC^{G,n}(X)
$$

\section{Orbifolds as Morita classes of groupoids}
We recall now the description of orbifolds as groupoids due to Moerdijk and Pronk \cite{MP}. Orbifolds were first introduced by Satake \cite{S} as a generalization of a manifold defined in terms of local quotients. The groupoid approach provides a global language to reformulate the notion of orbifold. This way of representing orbifolds allows for a natural generalization to the topological context, via topological groupoids.  We use the  same name orbifold for the topological version \footnote{Note that some authors call this notion \em{orbispace}.}.

 A topological groupoid $G_1\rightrightarrows G_0$ is {\it proper} if $(s,t):G_1\to G_0\times G_0$ is a proper map and it is a {\it foliation} groupoid if each isotropy group is discrete.

\begin{defn}
An \emph{orbifold groupoid} is a proper foliation groupoid.
\end{defn}

Given an orbifold groupoid $G_1\rightrightarrows G_0$, let $\orb(x)=t(s^{-1}(x))$ be the orbit through $x\in G_0$ and define the equivalence relation: $x\sim y$ if and only if $\orb(x)=\orb(y)$.
The orbit space ${G_0}/{\sim}$ is the space of all orbits with the quotient
topology. The orbit space of an orbifold groupoid is a locally compact Hausdorff space. Given an arbitrary locally compact Hausdorff space $X$ we can equip it  with an orbifold structure as follows.

\begin{defn} An \emph{orbifold structure} on a locally compact Hausdorff space $X$ is given by an orbifold groupoid $G_1\rightrightarrows G_0$ and a homeomorphism $h:{G_0}/{\sim}\to X$.
\end{defn}

If $\e$ is an essential equivalence between the groupoids
and $\tilde\e:{H_0}/{\sim}\to {G_0}/{\sim}$ is the induced homeomorphism between orbit spaces, we say that the composition $h\circ\tilde\e:{H_0}/{\sim}\to X$ defines an {\it equivalent} orbifold structure.

\begin{defn}  An \emph{orbifold}  $\X$
is a space $X$ equipped with a Morita equivalence class of orbifold
structures. A specific such structure, given by
$G_1\rightrightarrows G_0$ and  $h : {G_0}/{\sim} \to X $ is
a \emph{presentation} of the orbifold
$\X$.
\end{defn}

If two groupoids are Morita equivalent, then they define the same orbifold. Therefore any structure or invariant for orbifolds, if defined through groupoids, has to be invariant under Morita equivalence.

A large class of orbifolds, possibly all,  can be presented by a groupoid Morita equivalent to a translation groupoid $G\ltimes X$ with $G$ a compact group  acting  on $X$. Orbifolds that can be described this way are called representable. It is conjectured that all orbifolds are representable \cite{Henriques}.

From Corollary \ref{MoritaInvariant}, we have a well defined invariant for representable orbifolds:

\begin{defn}
Let $\X$ be a representable orbifold presented by the translation groupoid $G\ltimes X$ where $G$ is a compact Lie group and $X$ a metrizable space.
The {\em orbifold invariant topological complexity} of $\X$, ${\TC_{\cO}}(\X)$, is the invariant topological complexity of the groupoid $G\ltimes X$; that is ${\TC_{\cO}}(\X)={\TC}^G(X)$.
\end{defn}

\bibliographystyle{amsalpha}

\begin{thebibliography}{}

\bibitem{angelcolman} Andr\'es Angel and Hellen Colman, \emph{Equivariant topological complexities}, Topological Complexity and Related Topics, Contemporary Mathematics 702 (2018).

\bibitem{bayehsarkar} Marzieh Bayeh and Soumen  Sarkar, \emph{Higher equivariant and invariant topological complexity}, 	arXiv:1804.08006.

\bibitem{BlaszczykKaluba}
Zbigniew {B{\l}aszczyk} and Marek Kaluba, \emph{On equivariant and invariant topological complexity of smooth $\Z_p$-spheres}, Proc. Amer. Math. Soc. 145 (2017), 4075--4086.

\bibitem{bredon}  Glen Bredon, \emph{Introduction to Compact Transformation Groups}, Academic Press, New York, 1972

\bibitem{Clapp} M\'onica Clapp and Dieter Puppe, \emph{Invariants of the Lusternik-Schnirelmann type
and the topology of critical sets}, Trans. Amer. Math. Soc. 298 (1986) 603-620.

\bibitem{C}
Hellen Colman,
\emph{Equivariant LS-category for finite group actions},
Lusternik-Schnirelmann category and related topics (South Hadley, MA, 2001), 35-40,
Contemp. Math., 316, Amer. Math. Soc., Providence, RI, 2002.

\bibitem{CG}
Hellen Colman and Mark Grant, \emph{{Equivariant topological complexity}},
Algebr. Geom. Topol. \textbf{12} (2012), 2299--2316.


\bibitem{Clot03}
Octav Cornea, Gregory Lupton, John Oprea, and Daniel Tanr{\'e},
  \emph{Lusternik-{S}chnirelmann {C}ategory}, Mathematical Surveys and
  Monographs, vol. 103, American Mathematical Society, Providence, RI, 2003.

\bibitem{tomdieck}
Tammo tom Dieck,  \emph{ Transformation groups}, 1987.


\bibitem{Fa} Edward Fadell.\textit{ The equivariant Ljusternik-Schnirelmann method
for invariant functionals and relative cohomological index theories},
In
M\'ethodes Topologiques en Analyse Non-Lineaire, ed. A. Granas,
Montreal,
1985.


\bibitem{Far03} Michael Farber,
\emph{Topological complexity of motion planning}, Discrete Comput. Geom. 29. (2003), 211--221.

\bibitem{Henriques} Andre Henriques and David Metzler, \emph{Presentations of noneffective orbifolds}, Trans. Amer. Math. Soc. 356(2004), no. 6, 2481--2499.



\bibitem{wacek}
Wojciech Lubawski and Wac{\l}aw Marzantowicz, \emph{{Invariant topological
  complexity}}, Bull. London Math. Soc. \textbf{47} (2014), 101--117.

  \bibitem{Ma} Wac{\l}aw Marzantowicz, \textit{A G-Lusternik-Schnirelmann category of
    space with an action of a compact Lie group},  Topology  \textbf{ 28}
(1989),
403-412.

\bibitem{milnor}
John W.\ Milnor, \emph{Microbundles. I.},
Topology 3 (1964), suppl. 1, 53–-80.

\bibitem{MP}   Ieke Moerdijk \and Dorothea Pronk,
\emph{Orbifolds, Sheaves and Groupoids},  K-theory 13 (1997)
3--21.

\bibitem{MS}
Mitutaki Murayama and Kazuhisa Shimakawa,
\emph{Universal Equivariant Bundles}, Proc. AMS vol. 123 No. 4 (1995).

\bibitem{palais}Richard Palais, \emph{The classification of G-spaces}, Memoirs Amer. Math. Soc. 36 (1960).

\bibitem{pronk}
Dorette Pronk and Laura Scull, \emph{{Translation Groupoids and
Orbifold Cohomology}}, Canad. J. Math. Vol. 62 (3), 2010, 614--645.


\bibitem{rudyak} Yuli B. Rudyak. \emph{On higher analogs of topological complexity}, Topology Appl., 157(5):916--920, 2010.


\bibitem{S} Albert Schwarz,
 \emph{The genus of a fiber space},
 Amer.\ Math.\ Soc.\ Transl.  55 (1966), no. 2, 49--140.





\end{thebibliography}

\end{document}